\documentclass[preprint,times]{elsarticle}
\usepackage{mathtools,amssymb,amsfonts,ascmac}
\usepackage{theorem,multirow}
\usepackage{algorithm,algpseudocode}
\usepackage{url,xcolor,hhline}
\usepackage{graphicx}
\usepackage{mathptmx}
\SetSymbolFont{operators}{bold}{OT1}{txr}{bx}{n}
\SetSymbolFont{letters}{bold}{OML}{txmi}{bx}{it}
\SetSymbolFont{symbols}{bold}{OMS}{txsy}{bx}{n}
\SetSymbolFont{largesymbols}{bold}{OMX}{txex}{bx}{n}
\usepackage{bm}

\makeatletter
\renewcommand{\ALG@name}{Alg.}%
\makeatother

\journal{}

\begin{document}

\begin{frontmatter}

\title{On Implementation and Evaluation of Inverse Iteration Algorithm with compact WY Orthogonalization}
\author{Hiroyuki Ishigami}
\ead{hishigami@amp.i.kyoto-u.ac.jp}
\author{Kinji Kimura}
\ead{kkimur@amp.i.kyoto-u.ac.jp}
\author{Yoshimasa Nakamura}
\ead{ynaka@i.kyoto-u.ac.jp}

\address{Department of Applied Mathematics and Physics, Graduate School of Informatics, Kyoto University, Sakyo-ku, Kyoto 606 8501, Japan}


\begin{abstract}
A new inverse iteration algorithm that can be used to compute all the eigenvectors 
of a real symmetric tri-diagonal matrix on parallel computers is developed. 
The modified Gram-Schmidt orthogonalization is used 
in the classical inverse iteration. 
This algorithm is sequential and causes a bottleneck in parallel computing. 
In this paper, 
the use of the compact WY representation is proposed 
in the orthogonalization process of the inverse iteration with the Householder transformation. 
This change results in drastically reduced synchronization cost in parallel computing. 
The new algorithm is evaluated on both an 8-core and a 32-core parallel computer, 
and it is shown that 
the new algorithm is greatly faster than the classical inverse iteration algorithm 
in computing all the eigenvectors of matrices with several thousand dimensions. 
\end{abstract}


\end{frontmatter}

%
\section{Introduction}
The eigenvalue decomposition of a symmetric matrix, i.e., 
a decomposition into a product of matrices consisting of eigenvectors and eigenvalues, 
is one of the most important operations in linear algebra. 
It is used in vibrational analysis, image processing, data searches, etc. \par
Let us note that 
the eigenvalue decomposition of real symmetric matrices is reduced to 
that of real symmetric tri-diagonal matrices. 
Owing to recent improvements in the performance of computers 
equipped with multicore processors, 
we have had more opportunities to perform computation on parallel computers. 
As a result, 
there has been an increase in demand for an eigenvalue decomposition algorithm 
that can be effectively parallelized. \par
The inverse iteration algorithm is an algorithm 
for computing eigenvectors independently associated with mutually distinct eigenvalues. 
However, when we use this algorithm, we must reorthogonalize the eigenvectors 
if some eigenvalues are very close to each other. 
Adding this reorthogonalization process increases the computational cost. 
For this reorthogonalization, 
we have generally used the MGS (modified Gram-Schmidt) algorithm. 
However, this algorithm is sequential and inefficient for parallel computing. 
As a result, we are unable to maximize the performance of parallel computers. 
Hereinafter, 
we will refer to the inverse iteration algorithm with MGS as the classical inverse iteration. \par
We can also orthogonalize vectors by using the Householder transformation \cite{house}
and we call this precess the Householder orthogonalization algorithm. 
While the MGS algorithm is unstable in the sense 
that the orthogonality of the resulting vectors 
crucially depends on the condition number of the matrix \cite{yy}, 
the Householder algorithm is stable 
because its orthogonality does not depend on the condition number. 
The Householder algorithm is also sequential and ineffective for parallel computing, 
and its computational cost is higher than that of MGS. \par
In 1989, the Householder orthogonalization in terms of the compact WY representation 
was proposed  by R. Schreiber {\em et al} \cite{cwy}. 
By adopting this orthogonalization, 
stability and effective parallelization can be achieved. 
Hereafter, we refer to this algorithm as the compact WY orthogonalization algorithm. 
Yamamoto {\em et al.} \cite{yy} reformulated this algorithm 
for an incremental orthogonalization. 
Moreover, They showed that 
this algorithm achieves theoretically high accurate orthogonality 
and high scalability in parallel computing \cite{yy}. 
Here, the incremental orthogonalization is implemented on many numerical computation library. 
LAPACK(Linear Algebra PACKage) \cite{lapack} is one of the most popular libraries 
and all the code of LAPACK is implemented by using BLAS (Basic Linear Algebra Subroutines ) operations. 
The compact WY orthogonalization algorithm can be implemented by using BLAS. \par
In \cite{tinv}, 
authors have implemented the compact WY orthogonalization 
to the reorthogonalization process of inverse iteration for computing eigenvectors of a tri-diagonal matrix. 
It is shown \cite{tinv} that, in parallel computing, 
the new inverse iteration algorithm is faster than the classical one. \par
In this paper, 
we present two implementations: 
One is a new implementation of the compact WY orthogonalization algorithm based on BLAS. 
We focus on a mathematical structure of this algorithm 
and reformulate this algorithm. 
Therefore, using this new implementation, 
the computational cost of the compact WY orthogonalization can be reduced. 
The other is an implementation of the compact WY orthogonalization 
to the inverse iteration algorithm for a real symmetric tri-diagonal matrix. 
Thereafter, 
we perform the numerical experiments by computing all the eigenvectors using the second implementation 
and evaluate its performance. 
\section{Classical inverse iteration and its defect}
\subsection{Classical inverse iteration}
We consider the problem of computing eigenvectors of 
a real symmetric tri-diagonal matrix $T \in \mathbb{R}^{n \times n}$. 
Let $\lambda_j \in \mathbb{R}$ be eigenvalues of $T$ 
such that $\lambda_1 < \lambda_2 < \cdots < \lambda_n$. 
Let $\bm{v}_j \in \mathbb{R}^n$ be the eigenvector associated with $\lambda_j$. 
When $\tilde{\lambda_j}$, an approximate value of $\lambda_j$, 
and a starting vector $\bm{v}_j^{(0)}$ are given, 
we can compute an eigenvectors of $T$. 
To this end, we solve the following equation iteratively: 
\begin{equation}
  \left( T - \tilde{\lambda}_j I \right) \bm{v}_j^{(k)}  = \bm{v}_j^{(k-1)}. \label{eq:ii}
\end{equation}
Here $I$ is the $n$-dimensional identity matrix. 
If the eigenvalues of $T$ are mutually well-separated, 
$\bm{v}_j^{(k)}$, the solution of Eq.\eqref{eq:ii}, generically converges to the eigenvector associated with $\lambda_j$ 
as $k$ goes to $\infty$. 
The above iteration method is the inverse iteration. 
The computational cost of this method is of $O(mn)$ 
when we compute $m$ eigenvectors. 
In the implementation, we have to normalize the vectors $\bm{v}_j^{(k)}$ to avoid overflow. \par
When some of the eigenvalues are close to each other or there are clusters of eigenvalues of $T$, 
we have to reorthogonalize all the eigenvectors associated with such eigenvalues 
because they need to be orthogonal to each other. 
In the classical inverse iteration, 
we apply the MGS to this process and the computational cost of it is of $O(m^2 n)$. 
Therefore, 
when we compute eigenvectors of the matrix that has many clustered eigenvalues, 
the total computational cost increases significantly. 
In addition, 
the classical inverse iteration is implemented the Peters-Wilkinson method \cite{ii}. 
In this method, 
when the distance between the close eigenvalues is less than $10^{-3}\| T \|$, 
we regard them as members of the same cluster of eigenvalues, 
and we orthogonalize all of the eigenvectors associated with these eigenvalues. 
The classical inverse iteration algorithm is shown by Alg.\ref{alg:ii}, 
and $j_1$ denotes the index of the minimum eigenvalue of some cluster. 
This algorithm is implemented as DSTEIN 
in LAPACK \cite{lapack}. 
\begin{algorithm}[!t]
  \caption{Classical inverse iteration}
  \label{alg:ii}
  \begin{algorithmic}[1]
    \For{$j=1$ to $n$}
    \State Generate $\bm{v}_j^{(0)}$ from random numbers.
    \State $k=0$. 
    \Repeat
    \State $k \leftarrow k+1$.
    \State Normalize $\bm{v}_j^{(k-1)}$.
    \State Solve $\left( T - \tilde{\lambda}_j I \right) \bm{v}_j^{(k)}  = \bm{v}_j^{(k-1)}$ (Eq.\eqref{eq:ii}). 
    \If $|\tilde{\lambda_j} - \tilde{\lambda}_{j-1} | \le 10^{-3}\| T \|$, 
    \For{$i=j_1$ to $j-1$}
    \State $\bm{v}_j^{(k)} \leftarrow \bm{v}_j^{(k)}- \langle \bm{v}_j^{(k)}, \bm{v}_i \rangle \bm{v}_i$ 
    \EndFor
    \Else
    \State $j_1=j$.
    \EndIf
    \Until some condition is met.
    \State Normalize $\bm{v}_j^{(k)}$ to $\bm{v}_j$.
    \EndFor
  \end{algorithmic}
\end{algorithm}
\subsection{The defect of the classical inverse iteration}
The inverse iteration is a prominent method for computing eigenvectors, 
because we can compute eigenvectors independently. 
When there are many clusters in the distribution of eigenvalues, 
the inverse iteration can be parallelized by assigning each cluster to each core. \par
Let us consider the Peters-Wilkinson method in the classical inverse iteration. 
When the dimension of $T$ is greater than 1000, 
most of the eigenvalues are regarded as being in the same cluster \cite{mr}. 
In this case, 
we have to parallelize the inverse iteration 
with respect to not the cluster but the loop described from lines 2 to 16 in Alg.\ref{alg:ii}. 
This loop includes the iteration based on Eq.\eqref{eq:ii} and the orthogonalization of the eigenvectors. 
This orthogonalization process becomes 
a bottleneck of the classical inverse iteration 
with respect to the computational cost. 
The MGS algorithm is mainly based on a BLAS level-1 operation 
and it is a sequential algorithm. 
Because of this, 
when we compute all the eigenvectors on parallel computers, 
the number of synchronizations is of $O(m^2)$. 
Therefore, the MGS algorithm is ineffective in parallel computing. \par
In conclusion, 
the classical inverse iteration is an ineffective algorithm for parallel computing 
because the MGS algorithm is used in its orthogonalization process.
\section{Other orthogonalization algorithms}
In this section, 
we introduce alternative orthogonalization algorithms instead of the MGS algorithm. 
Now, we discuss the incremental orthogonalization of 
$\bm{v}_j \in \mathbb{R}^n$ to $\bm{q}_j \in \mathbb{R}^n$ 
($j=1$, $\dots$, $m$, $m \le n$). 
The incremental orthogonalization arises 
in the reorthogonalization process on the inverse iteration 
and it is defined as follows: 
$\bm{v}_j$ ($2 \le j \le m$) is not given in advance 
but is computed from $\bm{q}_1$, $\dots$, $\bm{q}_{j-1}$. \par
In the following, 
Let us define a vector $\bm{0}_{i}$ as the $i$-dimensional zero vector 
and matrices $V$, $Q \in \mathbb{R}^{n \times m}$ as 
$V = \left[ \bm{v}_1 \hspace{.5em} \cdots \hspace{.5em} \bm{v}_m \right]$, 
$Q = \left[ \bm{q}_1 \hspace{.5em} \cdots \hspace{.5em} \bm{q}_m \right]$.  
\subsection{Householder orthogonalization}
\begin{algorithm}[!t]
  \caption{Householder orthogonalization}
  \label{alg:house}
  \begin{algorithmic}[1]
    \For{$j=1$ to $m$}
    \State $\bm{u}_j  \leftarrow \left( I - t_1 \bm{y}_1 \bm{y}_1^\top \right) \bm{v}_j$
    \For{$i=2$ to $j-1$}
    \State $\bm{u}_j \leftarrow \left( I - t_i \bm{y}_i \bm{y}_i^\top \right) \bm{u}_j$
    \EndFor
    \State Compute $\bm{y}_j$ and $t_j$ by using $\bm{u}_j$ 
    \State $\bm{q}_j  \leftarrow \left( I - t_j \bm{y}_j \bm{y}_j^\top \right) \bm{e}_j$
    \For{$i=j-1$ to $1$}
    \State $\bm{q}_j \leftarrow \left( I - t_i \bm{y}_i \bm{y}_i^\top \right) \bm{q}_j$
    \EndFor
    \EndFor
  \end{algorithmic}
\end{algorithm}
The Householder orthogonalization, based on the Householder matrices, 
is one of the alternative orthogonalization methods. 
When vectors $\bm{u}_j$, 
$\bm{w}_j \in \mathbb{R}^n$ ($j=1$, $\dots$, $m$) 
satisfy $ \| \bm{u}_j \|_2 = \| \bm{w}_j \|_2$,  
there exists the orthogonal matrices $H_j$ called the Householder matrices satisfying 
$H_j H_j^\top = H_j^\top H_j = I$, $H_j \bm{u}_j = \bm{w}_j$ 
defined by $H_j = I -  t_j \bm{y}_j \bm{y}_j^\top$, 
$\bm{y}_j = \bm{u}_j - \bm{w}_j$, $t_j = 2 / \| \bm{y}_j \|_2^2$. 
The transformation from $\bm{u}_j$ to $\bm{v}_j$ by $H_j$ 
is called the Householder transformation. 
By using the Householder transformations. 
This orthogonalization algorithm is shown in Alg.\ref{alg:house}. 
The vector $\bm{y}_j$ is the vector 
in which the elements from 1 to $(j-1)$ are the same 
as the elements of $\bm{u}_j$ and the elements from $(j+1)$ to $n$ are zero. 
The vectors $\bm{u}_j$ and $\bm{w}_j$ are defined as follows: 
\begin{align*}
  \bm{u}_j 
  & = \begin{bmatrix} u_{1,j} & \cdots & u_{j-1,j} & u_{j,j} & u_{j+1,j} & \cdots & u_{n,j} \end{bmatrix}^\top \\
  & = H_{j-1} H_{j-2} \cdots H_2 H_1 \bm{v}_j , \\
  \bm{w}_j 
  & = \begin{bmatrix} u_{1,j} & \cdots & u_{j-1,j} & c_j & \bm{0}_{n-j}^\top \end{bmatrix}^\top, 
\end{align*}
where $u_{i,j}$ ($i=1$, $\dots$, $n$) is the $i$-th element of $\bm{u}_j$ and  
\begin{align*}
  c_j = -\operatorname{sgn} (u_{j,j})\sqrt{\sum_{i=j}^{n} u^2_{i,j}}. 
\end{align*}
Here, $\bm{y}_j$ and $t_j$ are computed as follows:
\begin{align}
  \bm{y}_j  
  = \bm{u}_j - \bm{w}_j 
  = \begin{bmatrix} \bm{0}_{j-1}^\top & u_{j,j} - c_j & u_{j+1,j} & \cdots & u_{n,j} \end{bmatrix}^\top, 
  t_j 
  = \frac{2}{\| \bm{y}_j \|_2^2}. 
  \label{eq:yt}
\end{align}
The vector $\bm{e}_j$ in Alg.\ref{alg:house} 
is the $j$-th vector of an $n$-dimensional identity matrix. \par
The orthogonality of the vectors $\bm{q}_j$ 
generated by the Householder orthogonalization 
does not depend on the condition number of $V$. 
Therefore, 
the Householder orthogonalization is more stable than MGS. 
On the other hand, 
being similar to MGS, 
it is a sequential algorithm,  that is mainly based on a BLAS level-1 operation. 
Its computational cost is about twice higher than that of MGS. 
Thus the Householder orthogonalization is an ineffective algorithm for parallel computing. 
\subsection{Compact WY orthogonalization}
In 1989, the Householder orthogonalization in terms of the compact WY representation 
was proposed by Schreiber and van Loan \cite{cwy}. 
Yamamoto and Hirota \cite{yy} reformulated this algorithm 
for the incremental orthogonalization. 
This study suggests that the Householder orthogonalization 
becomes capable of computation with a BLAS level-2 operation 
in terms of the compact WY representation. 
They also showed that 
this algorithm achieved theoretically high orthogonality 
and high scalability in parallel computing \cite{yy}. \par
Now, we consider the Householder orthogonalization in Alg.\ref{alg:house} 
and we introduce the compact WY representation. 
First, we define 
$Y_1 = [\bm{y}_1] \in \mathbb{R}^{n \times 1}$
and $T_1 = [ t_1 ] \in \mathbb{R}^{1 \times 1}$. 
Let us define matrices $Y_j \in \mathbb{R}^{n \times j}$ 
and upper triangular matrices $T_j \in \mathbb{R}^{j \times j}$ recursively as follows:
\begin{equation}
  Y_j =	
  \begin{bmatrix}
    Y_{j-1} & \bm{y}_j \\
  \end{bmatrix}, \hspace{.5em}
  T_j = 
  \begin{bmatrix}
    T_{j-1} & -t_j T_{j-1} Y_{j-1}^\top \bm{y}_j \\
    \bm{0}_{j-1}^\top  & t_j
  \end{bmatrix}.
  \label{eq:cwy1}				
\end{equation}
In this case, the following equation holds 
\begin{equation}
  H_1 H_2 \cdots H_j = I - Y_j T_j Y_j^\top.
  \label{eq:cwy2}
\end{equation}
As shown in Eq.\eqref{eq:cwy2}, 
we can rewrite the product of the Householder matrices $H_1 H_2 \cdots H_j$ 
in a simple block matrix form. 
Here $I - Y_j T_j Y_j^\top$ is called 
the compact WY representation of the product $H_1 H_2 \cdots H_j$ of the Householder matrices. 
Alg.\ref{alg:cwy} shows the compact WY orthogonalization algorithm. 
\begin{algorithm}[!t]
\caption{compact WY orthogonalization algorithm}
\label{alg:cwy}
\begin{algorithmic}[1]
  \State Compute $\bm{y}_1$ and $t_1$ by using $\bm{u}_1=\bm{v}_1$ 
  \State $Y_1 = \left[ \bm{y}_1 \right] $, $T_1 = \left[ t_1 \right]$ 
  \State $\bm{q}_1  \leftarrow \left( I - Y_1 T_1 Y_1^\top \right) \bm{e}_j$ 
  \For{$j=2$ to $m$}
  \State $\bm{u}_j \leftarrow \left( I - Y_{j-1} T_{j-1}^\top Y_{j-1}^\top \right) \bm{v}_j$ 
  \State Compute $\bm{y}_j$ and $t_j$ by using $\bm{u}_j$ 
  \State $Y_j = \begin{bmatrix} Y_{j-1} & \bm{y}_j \end{bmatrix}$, 
  $T_j = 
  \begin{bmatrix} 
    T_{j-1} & -t_j T_{j-1} Y_{j-1}^\top \bm{y}_j \\ 
    \bm{0}  & t_j 
  \end{bmatrix}$. 
  \State $\bm{q}_j  \leftarrow \left( I - Y_j T_j Y_j^\top \right) \bm{e}_j$ 
  \EndFor
  \end{algorithmic}
\end{algorithm}
\subsection{Implementation of compact WY orthogonalization}
In this subsection, 
we discuss the implementation of the compact WY orthogonalization algorithm 
using BLAS operations. 
In addition, 
we discuss a mathematical structure of this algorithm 
and present a new implementation of the compact WY orthogonalization 
for reducing the computational cost and the usage of memory. 
\subsubsection{Ordinary implementation of compact WY orthogonalization using BLAS}
Now we discuss the implementation of the compact WY orthogonalization 
based on line $5$ to $8$ in Alg.\ref{alg:cwy} using BLAS operations. 

For the adaptation of BLAS operations, 
we have to reformulate the formula of line $5$ as follows: 
\begin{align*}
  \bm{u}_j 
  & = \left( I - Y_{j-1} T_{j-1}^\top Y_{j-1}^\top \right) \bm{v}_j \\
  & = \bm{v}_j - Y_{j-1} T_{j-1}^\top Y_{j-1}^\top \bm{v}_j
\end{align*}
Now we can implement this formula by using BLAS as follows: 
\begin{align*}
  \begin{cases}
    \bm{u}_j \leftarrow \bm{v}_j & \text{(DCOPY)}\\
    \bm{v}_{j-1}^\prime \leftarrow Y_{j-1}^\top \bm{u}_j + 0 \cdot \bm{v}^\prime_{j-1} & \text{(DGEMV)}\\
    \bm{v}_{j-1}^\prime \leftarrow T_{j-1}^\top \bm{v}^\prime_{j-1} & \text{(DTRMV)}\\
    \bm{u}_j \leftarrow (-1) \cdot Y_{j-1} \bm{v}_{j-1}^\prime + \bm{u}_j & \text{(DGEMV)}
  \end{cases},
\end{align*} 
where $\bm{v}^\prime_{j-1} \in \mathbb{R}^{j-1}$. 
We set the initial address of $\bm{v}^\prime_{j-1}$ assigned on CPU memory 
to correspond to that of $\bm{v}_j$. 
DCOPY denotes the copying operation of a vector $\bm{x}$ to a vector $\bm{y}$: 
$\bm{y} \leftarrow \bm{x}$. 
DGEMV means the matrix-vector operation: 
$\bm{y} \leftarrow \alpha A\bm{x} + \beta \bm{y}$, 
where $A$ is a general rectangular matrix. 
DTRMV denotes the matrix-vector product: 
$\bm{x} \leftarrow T\bm{x}$, 
where $T$ is a triangular matrix. 

Next, on line 6,
we compute $\bm{y}_j$ and $t_j$ based on Eq.\eqref{eq:yt}. 
These computations is mainly performed by using BLAS level-1 operations 
and its computational cost is relatively lower. 
we implement the computation of $\bm{y}_j$ and $t_j$ as follows: 
\begin{align*}
  \begin{cases}
    y_{i,j} \leftarrow 0, \text{ ($i=1$, $\dots$, $j-1$)} & \\
    y_{i,j} \leftarrow u_{i,j}, \text{ ($i=j$, $\dots$, $n$)} & \text{(DCOPY)}\\
    y_{j,j} \leftarrow u_{j,j} - c_j, \hspace{.5em} 
    c_j = -\operatorname{sgn} (u_{j,j})\sqrt{\sum_{i=j}^{n} u^2_{i,j}} 
    & \text{(DNRM2)} \\
    t_j \leftarrow 2/\| \bm{y}_j \|^2_2 & \text{(DNRM2)}
  \end{cases},
\end{align*}
where $y_{i,j}$ ($i=1$, \dots, $n$) is the $i$-th column element of $\bm{y}_j$. 
DNRM2 denotes the computation of the $2$-norm of a vector. 

On line 7, 
updating $Y_j$ and $t_j$ can be done easily. 
Now, let $\bm{\hat{t}}_j \in \mathbb{R}^{j-1}$ 
be $\bm{\hat{t}}_j = -t_j T_{j-1} Y_{j-1}^\top \bm{y}_j$.  
Note that $\bm{\hat{t}}_j$ is implemented by using BLAS as follows: 
\begin{align*}
  \begin{cases}
    \bm{\hat{t}}_j \leftarrow (-t_j) Y_{j-1}^\top \bm{y}_j + 0 \cdot \bm{\hat{t}}_j & \text{(DGEMV)}\\ 
    \bm{\hat{t}}_j \leftarrow T_{j-1} \bm{\hat{t}}_j & \text{(DTRMV)}
  \end{cases}. 
\end{align*}

At last, on line 8, 
we can reformulate as follows: 
\begin{align*}
  \bm{q}_j 
  & = \left( I - Y_j T_j Y_j^\top \right) \bm{e}_j \\
  & = \bm{e}_j - Y_j T_j Y_j^\top \bm{e}_j.  
\end{align*}
Here, 
the matrix-vector product $Y_j^\top \bm{e}_j$ can be simplified as follows: 
\begin{align*}
  Y_j^\top \bm{e}_j 
  = \begin{bmatrix} y_{j,1} \\ \vdots \\ y_{j,j} \end{bmatrix}  
\end{align*}. 
This computation can be performed only by copying the $j$-th column of $Y_j$ to some vector. 
Therefore we can implement the formula of line 8 using BLAS as follows:
\begin{align*}
  \begin{cases}
    \bm{q}_j \leftarrow \bm{e}_j & \text{(DCOPY)} \\
    \bm{v}_j^\prime \leftarrow \begin{bmatrix} y_{j,1} & \cdots & y_{j,j} \end{bmatrix} & \text{(DCOPY)}\\
    \bm{v}_j^\prime \leftarrow T_j^\top \bm{v}_j^\prime & \text{(DTRMV)}\\
    \bm{q}_j \leftarrow (-1) \cdot Y_j \bm{v}_j^\prime + \bm{q}_j & \text{(DGEMV)}
  \end{cases},
\end{align*}
where $\bm{v}^\prime_j \in \mathbb{R}^j$, $\bm{q}_j \in \mathbb{R}^n$. 
We set the initial address of $\bm{v}^\prime_j$, $\bm{q}_j$ assigned on CPU memory 
to correspond to that of $\bm{u}_j$, $\bm{v}_j$, respectively. 

The computational cost of the above compact WY orthogonalization algorithm is almost $4m^2n+m^3$. 
In the worst case, i.e., $m=n$, the computational cost is $5n^3$. 

In addition, 
for this implementation, 
we have to use almost $mn + m^2$ CPU memory 
because $Y_m$ use $mn$ and $T_m$ use $m^2$ domain. 
\subsubsection{New implementation of compact WY orthogonalization using BLAS}

In the above section, 
we discuss the ordinary implementation of the compact WY orthogonalization algorithm. 
Now we focus on the mathematical structure of this algorithm 
and present the new implementation of the compact WY orthogonalization 
which has the less computational cost than the ordinary one has. 

Before the formula of line 5 in Alg.\ref{alg:cwy}, 
let us consider the formula of line 6. 
From Eq.\eqref{eq:yt}, 
we can strictly compute $t_j$ as follows: 
Since 
\begin{align*}
  c_j & = -\operatorname{sgn} \left( u_{j,j} \right) \sqrt{\sum_{i=j}^{n} u^2_{i,j}},  
\end{align*}
we have 
\begin{align*}
  \| \bm{y}_j \|_2^2 
  & = \left( u_{j,j} - c_j \right)^2 + \sum_{i=j+1}^{n} u^2_{i,j} \\
  & = \sum_{i=j}^{n} u^2_{i,j} - 2u_{j,j} c_j + c_j^2 \\
  & = 2 ( c_j^2 - u_{j,j} c_j ). 
\end{align*}
Hence, we have 
\begin{align*}
  t_j & = \frac{2}{\| \bm{y}_j \|_2^2} = \frac{1}{c_j^2 - u_{j,j} c_j}. 
\end{align*}
From this fact and the definition of $\bm{y}_j$ and $c_j$, 
we need not compute the elements from $1$ to $(j-1)$ of $\bm{u}_j$ in actual.   
Therefore we compute only the elements from $j$ to $n$ of $\bm{u}_j$ 
so that the formula of line 5 is reduced as follows: 
\begin{align*}
  \bm{\hat{u}}_j = \bm{\hat{u}}_j - \hat{Y}_{j-1} T_{j-1}^\top Y_{j-1}^\top \bm{v}_j, 
\end{align*}
where $\bm{\hat{u}}_j \in \mathbb{R}^{n-(j-1)}$ is 
$\bm{\hat{u}}_j = \begin{bmatrix} u_{j,j} & \cdots & u_{n,j} \end{bmatrix}^\top$. 

Here, 
we focus on the structure of $\bm{y}_j$. 
From Eq.\eqref{eq:yt}, 
$\bm{y}_j$ ($j=2$, $\dots$, $m$) can be represented as the block vector of the form:  
\begin{align*}
  \bm{y}_j = 
  \begin{bmatrix}
    \bm{0}_{j-1} \\
    \bm{\hat{y}}_{j} 
  \end{bmatrix},  
\end{align*}
where $\bm{\hat{y}}_{j} \in \mathbb{R}^{n-(j-1)}$ is the vector of nonzero elements of $\bm{y}_j$. 
From this fact, 
$Y_j$ can be represented as the following block matrix:  
\begin{align*}
  Y_j = 
  \begin{bmatrix}
    L_j \\ 
    \hat{Y}_j
  \end{bmatrix}, 
\end{align*}
where $L_j \in \mathbb{R}^{j \times j}$ is a lower triangular matrix 
and $\hat{Y}_j \in \mathbb{R}^{(n-j) \times j}$ is generally a dense rectangular matrix. 
In addition, 
let us consider $\bm{v}_j$ as the block vector of the form: 
\begin{align*}
  \bm{v}_j 
  & = 
  \begin{bmatrix}
    \bm{\check{v}}_j \\
    \bm{\hat{v}}_j 
  \end{bmatrix}, 
\end{align*}
where $\bm{\check{v}}_j \in \mathbb{R}^{j-1}$, $\bm{\hat{v}}_j \in \mathbb{R}^{n-(j-1)}$. 

By using these block form of $\bm{v}_j$ and $Y_j$, 
we can reduce the computational cost of the matrix-vector product $Y_{j-1}^\top \bm{v}_j$ through 
\begin{align*}
  Y_{j-1}^\top \bm{v}_j = 
  \begin{bmatrix}
    L_{j-1} \\ 
    \hat{Y}_{j-1}
  \end{bmatrix}^\top 
  \begin{bmatrix}
    \bm{\check{v}}_j \\
    \bm{\hat{v}}_j 
  \end{bmatrix}
  = L_{j-1}^\top \bm{\check{v}}_j + \hat{Y}_{j-1}^\top \bm{\hat{v}}_j. 
\end{align*}
Therefore, the formula of $\bm{\hat{u}}_j$ can be simplified as follows: 
\begin{align*}
  \bm{\hat{u}}_j & = \bm{\hat{u}}_j - \hat{Y}_{j-1} T_{j-1}^\top 
  \left( L_{j-1}^\top \bm{\check{v}}_j + \hat{Y}_{j-1}^\top \bm{\hat{v}}_j \right).  
\end{align*}
This formula can be implemented by using BLAS as follows: 
\begin{align*}
  \begin{cases}
    \bm{\hat{u}}_j \leftarrow \bm{\hat{v}}_j & \text{(DCOPY)}\\
    \bm{\check{v}}_j \leftarrow L_{j-1}^\top \bm{\check{v}}_j & \text{(DTRMV)}\\
    \bm{\check{v}}_j \leftarrow \hat{Y}_{j-1}^\top \bm{\hat{v}}_j + \bm{\check{v}}_j & \text{(DGEMV)}\\
    \bm{\check{v}}_j \leftarrow T_{j-1}^\top \bm{\check{v}}_j & \text{(DTRMV)}\\
    \bm{\hat{u}}_j \leftarrow (-1) \cdot \hat{Y}_{j-1} \bm{\check{v}}_j + \bm{\hat{u}}_j & \text{(DGEMV)}
  \end{cases}.
\end{align*}

From the above discussion, 
the computation on line 6 is implemented by using BLAS as follows: 
\begin{align*}
  \begin{cases}
    y_{i,j} \leftarrow u_{i,j}, \text{ ($i=j$, $\dots$, $n$)} & \text{(DCOPY)}\\
    y_{j,j} \leftarrow u_{j,j} - c_j, \hspace{.5em} 
    c_j = -\operatorname{sgn} (u_{j,j})\sqrt{\sum_{i=j}^{n} u^2_{i,j}} 
    & \text{(DNRM2)} \\
    t_j \leftarrow 1/\left(c_j^2 - u_{j,j} c_j \right) & 
  \end{cases}.
\end{align*}

On line 7, 
we can also reduce the computational cost of $\bm{\hat{t}}_j$ through 
\begin{align*}
  \bm{\hat{t}}_j 
  & = -t_j T_{j-1} Y_{j-1}^\top \bm{y}_j \\
  & = -t_j T_{j-1} 
  \begin{bmatrix}
    L_{j-1} \\ 
    \hat{Y}_{j-1}
  \end{bmatrix}^\top 
  \begin{bmatrix}
    \bm{0}_{j-1} \\
    \bm{\hat{y}}_{j} 
  \end{bmatrix} \\
  & = -t_j T_{j-1} 
  \left( L_{j-1}^\top \bm{0}_{j-1} + \hat{Y}_{j-1}^\top \bm{\hat{y}}_j \right) \\
  & = -t_j T_{j-1} \hat{Y}_{j-1}^\top \bm{\hat{y}}_j.  
\end{align*}
This formula can be implemented by using BLAS as follows: 
\begin{align*}
  \begin{cases}
    \bm{\hat{t}}_j \leftarrow (-t_j) \hat{Y}_{j-1}^\top \bm{\hat{y}}_j + 0 \cdot \bm{\hat{t}}_j & \text{(DGEMV)}\\ 
    \bm{\hat{t}}_j \leftarrow T_{j-1} \bm{\hat{t}}_j & \text{(DTRMV)}
  \end{cases}.
\end{align*}

At last, on line 8, 
even if the sign of the orthogonal vector $\bm{q}_j$ is reversed, 
the orthogonality along with other vectors is not changed. 
Therefore, 
we can reformulate $\bm{q}_j$ as $\bm{q}_j = \left( Y_j T_j Y_j^\top - I \right) \bm{e}_j$. 
In addition, 
let us consider $\bm{q}_j$ as the following block vector:
\begin{align*}
  \bm{q}_j 
  & = 
  \begin{bmatrix}
    \bm{\check{q}}_j \\
    \bm{\hat{q}}_j 
  \end{bmatrix}, 
\end{align*}
where $\bm{\check{q}}_j \in \mathbb{R}^j$, $\bm{\hat{q}}_j \in \mathbb{R}^{n-j}$. 
These are reformulated as follows: 
\begin{align*}
  \begin{bmatrix}
    \bm{\check{q}}_j \\
    \bm{\hat{q}}_j 
  \end{bmatrix} 
  & = 
  \begin{bmatrix}
    L_j T_j Y_j^\top \bm{e}_j \\
    \hat{Y}_j T_j Y_j^\top \bm{e}_j 
  \end{bmatrix} - 
  \begin{bmatrix}
    \bm{\check{e}}_j \\
    \bm{0}_{n-j}
  \end{bmatrix}, 
\end{align*}
where $\bm{\check{e}}_j$ is the $j$-th vector of the $j$-dimensional identity matrix. 
Therefore this formula can be implemented by using BLAS as follows: 
\begin{align*}
  \begin{cases}
    \bm{x}_j \leftarrow \begin{bmatrix} y_{j,1} & \cdots & y_{j,j} \end{bmatrix} & \text{(DCOPY)}\\
    \bm{x}_j \leftarrow T_j^\top \bm{x}_j & \text{(DTRMV)}\\
    \bm{\check{q}}_j \leftarrow \bm{x}_j & \text{(DCOPY)}\\
    \bm{\check{q}}_j \leftarrow L_j \bm{\check{q}}_j & \text{(DTRMV)}\\
    \bm{\hat{q}}_j \leftarrow \hat{Y}_j \bm{x}_j + 0 \cdot \bm{\hat{q}}_j & \text{(DGEMV)}\\
    q_{j,j} \leftarrow q_{j,j} - 1 &
  \end{cases},
\end{align*}
where $\bm{x}_j \in \mathbb{R}^j$ is assigned on workspace memory. 

When the above implementation is adapted, 
the highest order of the computational cost of the compact WY algorithm reduced to $4m^2n-m^3$. 
In the worst case, i.e., $m=n$,
the computational cost of the new implementation of the compact WY algorithm is almost $3 n^3$. 

\begin{figure}[!t]
  \begin{center}
    \includegraphics[scale=.5]{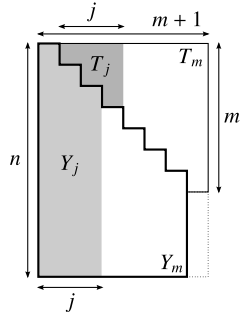}
  \end{center}
  \caption{Assignment model for $Y_j$ and $T_j$}
  \label{figure:assignment}
\end{figure}

In addition, 
our implementation have not to be referred any zero elements of $Y_j$ and $T_j$. 
Therefore, if $Y_j$ and $T_j$ are assigned on a CPU memory like Alg.\ref{figure:assignment}, 
the use of memory can be reduced to almost $n(m+1)$, 
\subsection{Comparison of the orthogonalization algorithms}
The compact WY orthogonalization has a stable orthogonality 
arising from the Householder transformations, 
and its numerical computation is mainly performed by BLAS level-2 operations. 
As a result, 
this orthogonalization has a better stability and a sophisticated orthogonality, 
and it is more effective for parallel computing than MGS. 
Table \ref{table:perform} displays the differences 
in performance of the orthogonalization methods mentioned above. 
In this table, 
{\em Computation} denotes the order of the computational cost. 
{\em Synchronization} means the order of the number of synchronizations. 
{\em Orthogonality} indicates the norm  $\| Q^\top Q - I \|$ 
and $\epsilon$ denotes the machine epsilon 
and $\kappa(V)$ is the condition number of $V$. 
\begin{table}[!t]
\centering
\footnotesize
\caption{Comparison of the orthogonalization methods \cite{perform} \cite{yy}}
\label{table:perform}
\begin{tabular}{rccc}\hline
      orthogonalization & {\em Computation} & {\em Synchronization} & {\em Orthogonality} \\ \hline
      MGS               & $2m^2 n$   & $O(m^2)$              & $O(\epsilon \kappa(V) )$ \\ 
      Householder       & $4m^2 n$   & $O(m^2)$              & $O(\epsilon)$ \\
      compact WY        & $4m^2 n + m^3$ & $O(m)$                & $O(\epsilon)$ \\
      new compact WY    & $4m^2 n - m^3$ & $O(m)$                & $O(\epsilon)$ \\\hline
\end{tabular}
\end{table}
\section{Inverse iteration algorithm with compact WY orthogonalization}
Authors have proposed an alternative inverse iteration algorithm in \cite{tinv}. 
This algorithm is based on the classical inverse iteration algorithm implemented in DSTEIN 
and we change the orthogonalization process of it from MGS to the compact WY orthogonalization 
that is described on Sec. 3.3.1.
In addition, 
it is shown that this algorithm is faster than the classical inverse iteration one 
in parallel computing \cite{tinv}. \par
Now we present an even faster inverse iteration algorithm with the compact WY orthogonalization. 
This compact WY orthogonalization is implemented on the way of Sec. 3.3.2. 
The new algorithm is described in Alg.\ref{alg:cwyii}.
Let us name the new code DSTEIN-cWY. 
\begin{algorithm}[!t]
  \caption{compact WY inverse iteration}
  \label{alg:cwyii}
  \begin{algorithmic}[1]
    \For{$j=1$ to $n$}
    \State Generate $\bm{v}_j^{(0)}$ from random numbers.
    \State $k=0$
    \Repeat
    \State $k \leftarrow k+1$.
    \State Normalize $\bm{v}_j^{(k-1)}$.
    \State Solve $\left( T - \tilde{\lambda}_j I \right) \bm{v}_j^{(k)}  = \bm{v}_j^{(k-1)}$.
    \If{$|\tilde{\lambda_j} - \tilde{\lambda}_{j-1} | \le 10^{-3}\| T \|$, }
    \State $j_c \leftarrow j - j_1$. 
    \If{$j_c=1$ and $k=1$, }
    \State Compute $Y_1 = [ \bm{y}_1 ]$ and $T_1 = [ t_1 ]$ by using $\bm{v}_{j_1}$.
    \EndIf
    \State $\bm{u}_{j_c+1} = \left( I - Y_{j_c} T_{j_c}^\top Y_{j_c}^\top \right) \bm{v}_j^{(k)}$.
    \State Compute $\bm{y}_{j_c+1}$ and $t_{j_c+1}$ by using $\bm{u}_{j_c+1}$.
    \State $Y_{j_c+1} = 
    \begin{bmatrix} 
      Y_{j_c} & \bm{y}_{j_c+1} 
    \end{bmatrix}$, 
    $T_{j_c+1} = 
    \begin{bmatrix} 
      T_{j_c} & -t_{j_c+1} T_{j_c} Y_{j_c}^\top \bm{y}_{j_c+1} \\ 
      \bm{0}_{j_c}^\top  & t_{j_c+1} 
    \end{bmatrix}$.
    \State $\bm{v}_j^{(k)} \leftarrow \left( I - Y_{j_c+1} T_{j_c+1} Y_{j_c+1}^\top \right) \bm{e}_{j_c+1}$.
    \Else
    \State $j_1 \leftarrow j$.
    \EndIf
    \Until{Some condition is met.}
    \State Normalize $\bm{v}_j^{(k)}$ to $\bm{v}_j$.
    \EndFor
  \end{algorithmic}
\end{algorithm}

Next, 
we explain an application of the new implementation of the compact WY orthogonalization 
to the inverse iteration. 
Differences between DSTEIN-cWY and DSTEIN is as follow: 
For the classical inverse iteration algorithm, 
we need not know the index $j_c$ which denotes the $j_c$-th eigenvalue of the cluster 
in computing the eigenvector associated with it. 
However, 
we must know the index for the compact WY orthogonalization 
when we compute and update $T_j$, $Y_j$. 
To overcome the above difficulty, 
we introduce a variable $j_c$ on line 9, and we can recognize it. 
This introduction of $j_c$ enables us to execute the intended program. \par
In the classical inverse iteration algorithm, 
we need not know the first eigenvalue $\lambda_{j_1}$ of the cluster. 
However, we must compute $\bm{y}_1$ and $t_1$ 
in the new inverse iteration algorithm. 
Therefore, 
at the starting point of the computation of the eigenvector 
associated with the second eigenvalue $\lambda_{j_1 + 1}$, 
we compute $T_1 = [t_1]$, $Y_1 = [\bm{y}_1]$ by using $\bm{v}_{j_1}$. 
At this time, 
because $\bm{v}_{j_1}$ is a normalized vector 
so that it equals to $( I - Y_1 T_1 Y_1^\top ) \bm{e}_1$, 
we need not compute $\bm{v}_{j_1}$ it again. 
\section{Numerical experiments}
We describe some numerical experiments performed 
by using DSTEIN and DSTEIN-cWY on parallel computers, 
and we compare the computation time. 
Here DSTEIN of LAPACK is based on the classical inverse iteration, 
and DSTEIN-cWY makes use of the new inverse iteration presented in the previous section. 
\subsection{Contents of the numerical experiments}
We report computations of all the eigenvectors associated 
with eigenvalues of some matrices by using DSTEIN and DSTEIN-cWY on parallel computers, 
and we compare the elapsed time. 
In these experiments, 
we compute the approximate eigenvalues by using LAPACK's program DSTEBZ, 
which is capable of computing eigenvalues using the bisection method. 
We record the elapsed time for DSTEIN and DSTEIN-cWY using SYSTEM\_CLOCK, 
which is the internal function of Fortran. \par
In the experiments, 
we use two computers equipped with multicore CPUs, 
and we implement those algorithms by using GotoBLAS2 \cite{goto}, 
which is implemented to parallelize BLAS operations by assigning them to each CPU core. 
Table \ref{table:2} shows the specifications of two computers. 
\begin{table}[!t]
\centering
\caption{The specification of Computer 1 and 2}
\label{table:2}
\small
\begin{tabular}{rcc}\hline
               &    Computer 1      &  Computer 2         \\ \hline
               & AMD Opteron 2.0GHz &  Intel Xeon 2.93GHz \\
      \raisebox{.5\normalbaselineskip}[0pt][0pt]{CPU}
               & 32cores(8cores$\times$4)&  8cores(4cores$\times$2)  \\
      RAM      & 256GB              &  32GB               \\
      Compiler & Gfortran-4.4.5     &  Gfortran-4.4.5     \\
      LAPACK   & LAPACK-3.3.0       &  LAPACK-3.3.0       \\
      BLAS     & GotoBLAS2-1.13     & GotoBLAS2-1.13      \\ \hline
\end{tabular}
\end{table}
As experimental matrices, we use symmetric tri-diagonal matrices of three types. 
Type 1 is a tri-diagonal random matrix, of which elements are set to the random number of $[0,1)$. 
It is shown that 
the eigenvalues of a tri-diagonal random matrix are divided into a few clusters in the sense of Peters-Wilkinson method\cite{ii}. 
and most of eigenvalues are included in the biggest one of the clusters  
if the dimension $n$ of a random matrix becomes larger. 
The tri-diagonal matrix of Type 2 is defined as follows: 
\begin{equation}
  T	=
  \begin{bmatrix}
    1 & 1 &        & \\
    1 & 1 & 1      & \\
      & 1 & \ddots & \ddots \\
      &   & \ddots & \ddots  & 1 \\
      &   &        & 1       & 1 
  \end{bmatrix}.
\end{equation}
All the eigenvalues of Type 2 matrix with large dimensions 
are included in the same cluster in the sense of Peters-Wilkinson method. 
Type 3 is the glued-Wilkinson matrices $W_g^\dagger$. 
$W_g^\dagger$ consists of the block matrix $W_{21}^\dagger \in \mathbb{R}^{21 \times 21}$ 
and the scalar parameter $\delta \in \mathbb{R}$ and is defined as follow: 
\begin{equation}
  W_g^\dagger =
  \left[
    \begin{array}{cc|cc|cc|cc}
      \multicolumn{2}{c|}{\raisebox{-0.5em}[0pt][0pt]{\normalsize $W_{21}^\dagger$}}& & & & & & \\
      & & \delta & & & & & \\\hline
      & \delta & \multicolumn{2}{c|}{\raisebox{-0.5em}[0pt][0pt]{\normalsize $W_{21}^\dagger$}} & & & & \\
      & & & & \delta & & & \\\hline
      & & & \delta & \ddots & \ddots & & \\
      & & & & \ddots & \ddots & \delta & \\\hline
      & & & & & \delta & \multicolumn{2}{c}{\raisebox{-0.5em}[0pt][0pt]{\normalsize $W_{21}^\dagger$}} \\
      & & & & & & &
    \end{array}
    \right], 
\end{equation}
where $W_{21}^\dagger$ is defined by 
\begin{equation}
  W_{21}^\dagger 
  =
  \begin{bmatrix}
    10 & 1 &        &        &         &	\\
    1  & 9 & 1      &        &         &	\\
    & 1 & \ddots & \ddots &         & 	\\
    &   & \ddots & 0      & \ddots  &  	\\
    &   &        & \ddots & \ddots  & 1	\\
    &   &        &        & 1       & 10
  \end{bmatrix}, 
\end{equation}
and $\delta$ satisfies $0 < \delta < 1$ and is also the semi-diagonal element of $W_g^\dagger$. 
Since $W_g^\dagger$ is real symmetric tri-diagonal and its semi-diagonal elements are nonzero, 
all the eigenvalues of $W_g^\dagger$ are real 
and they are divided into 21 clusters of close eigenvalues. 
When $\delta$ is small, 
the distance between the minimum and maximum eigenvalues in any cluster is small. 
In our experiments, we set $\delta = 10^{-4}$. 
Computing eigenvalues and eigenvectors of the glued-Wilkinson matrix 
is one of the benchmark problems of eigenvalue decomposition. 
For example, 
the glued-Wilkinson matrix was used to evaluate 
the performance of matrix eigenvalue algorithms \cite{gw1} \cite{gw2}. \par
\subsection{Results of the experiments}
\par
\begin{figure*}[!t]
 \begin{minipage}{0.5\hsize}
  \begin{center}
   \includegraphics[scale=.5]{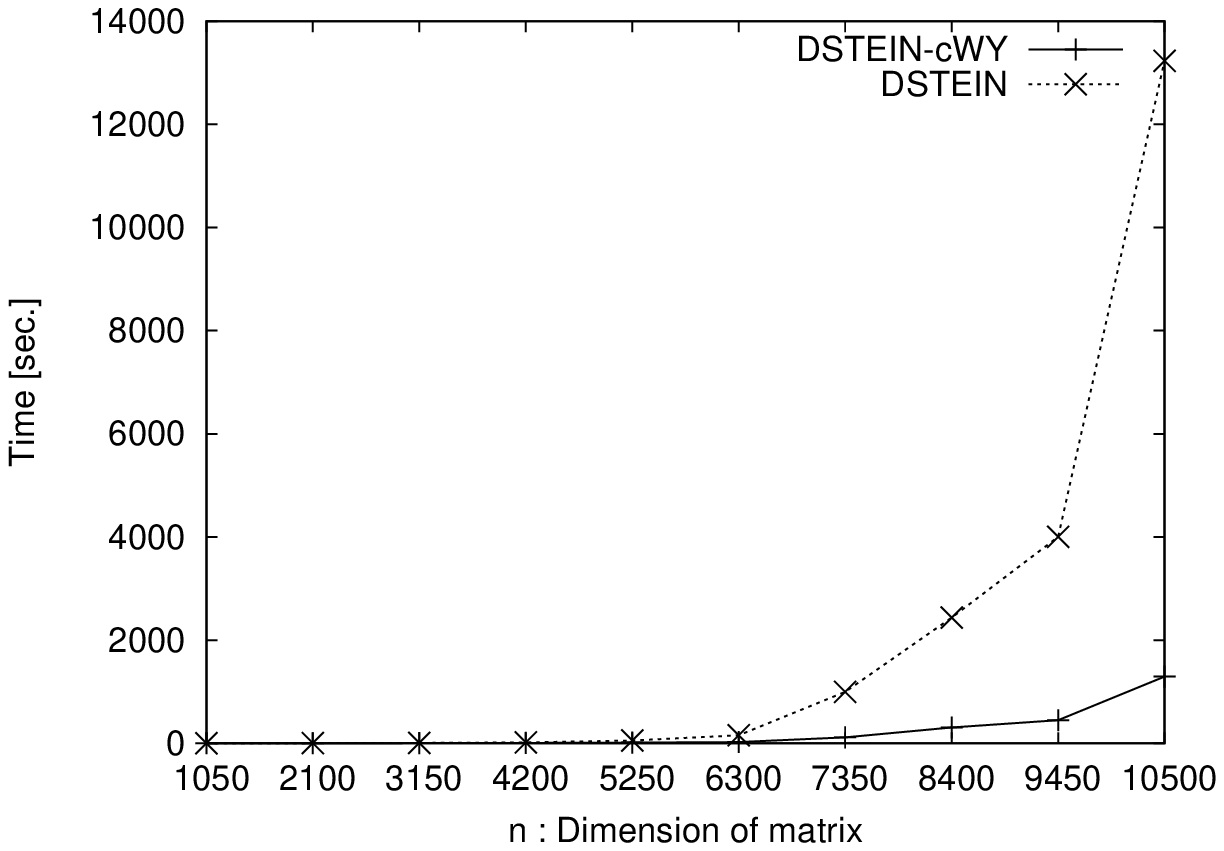}
  \end{center}
 \end{minipage}
 \begin{minipage}{0.5\hsize}
  \begin{center}
   \includegraphics[scale=.5]{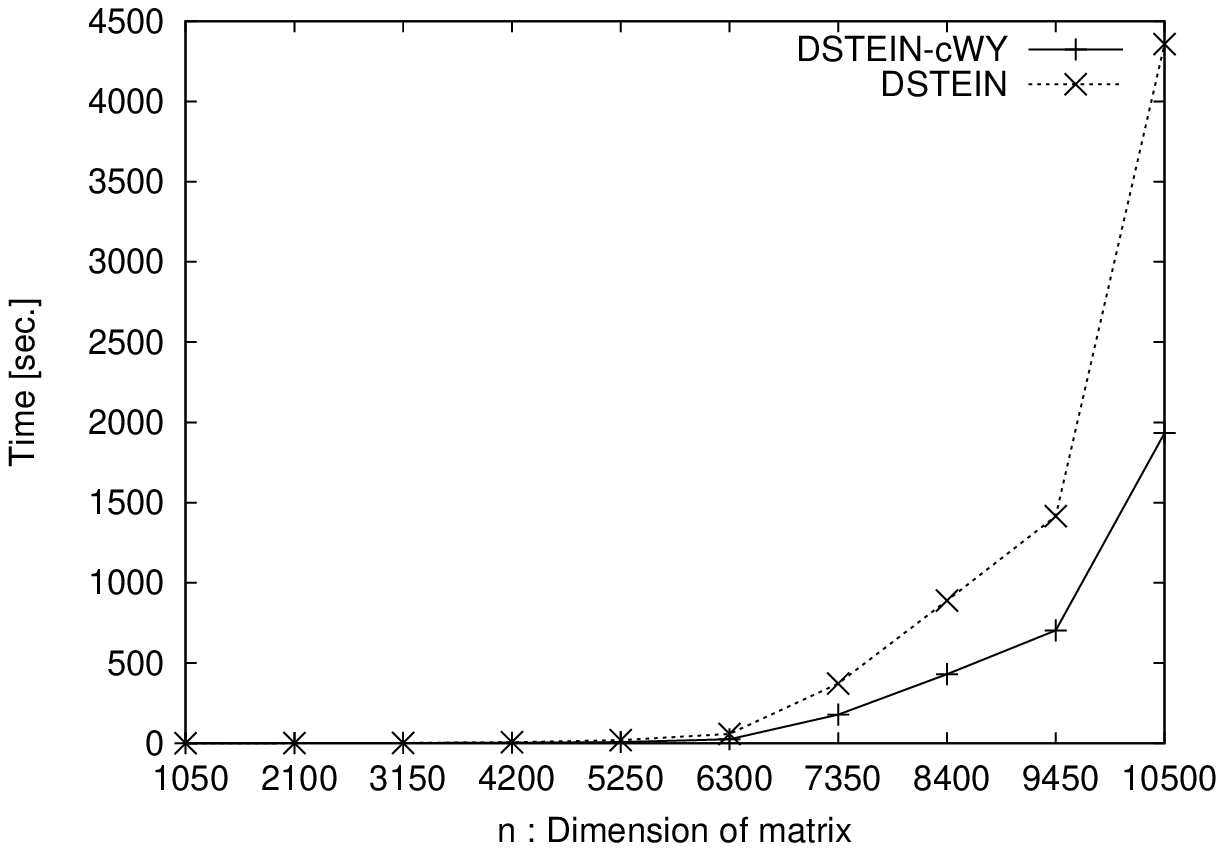}
  \end{center}
 \end{minipage}
  \caption{Dimension $n$ of Type 1 matrix and the computation time by DSTEIN and DSTEIN-cWY. 
the left graph corresponds to Computer 1 and the right Computer 2.}
  \label{graph:1}
 \begin{minipage}{0.5\hsize}
  \begin{center}
   \includegraphics[scale=.5]{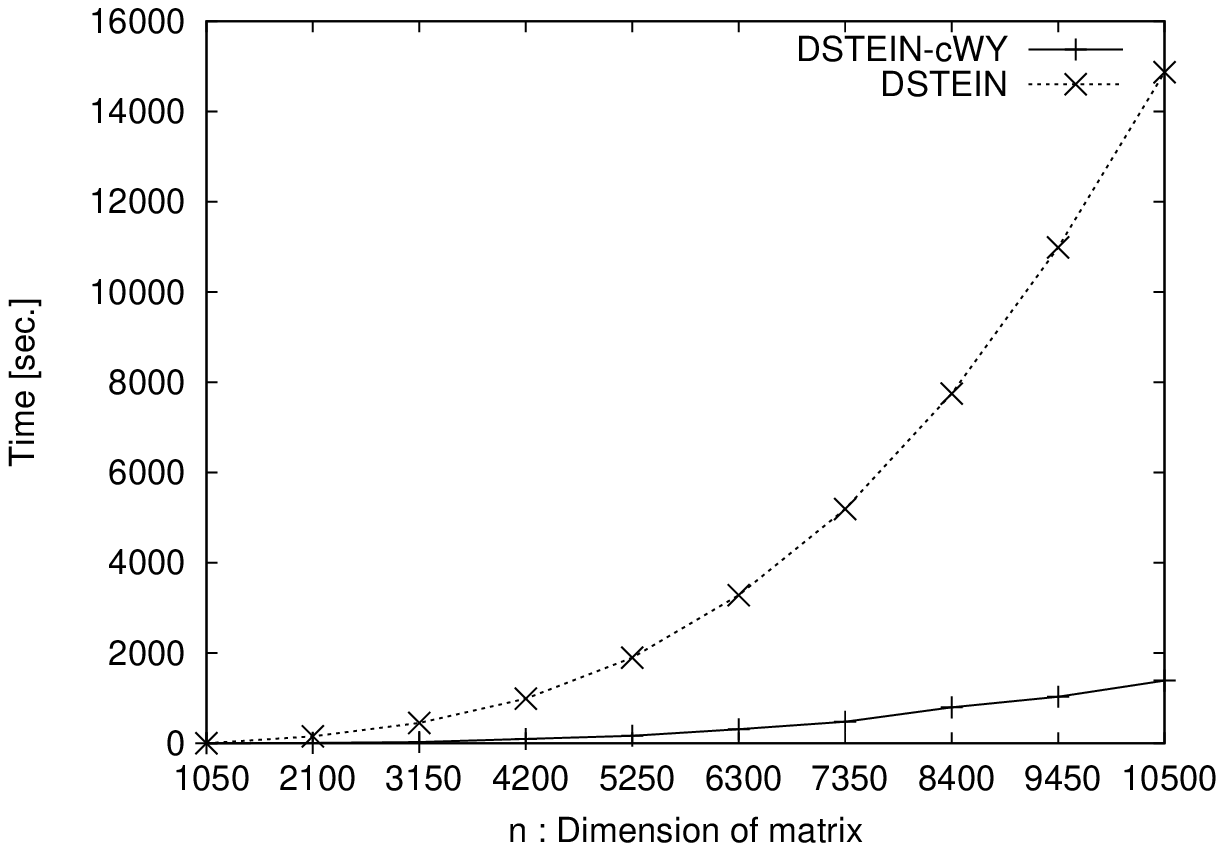}
  \end{center}
 \end{minipage}
 \begin{minipage}{0.5\hsize}
  \begin{center}
   \includegraphics[scale=.5]{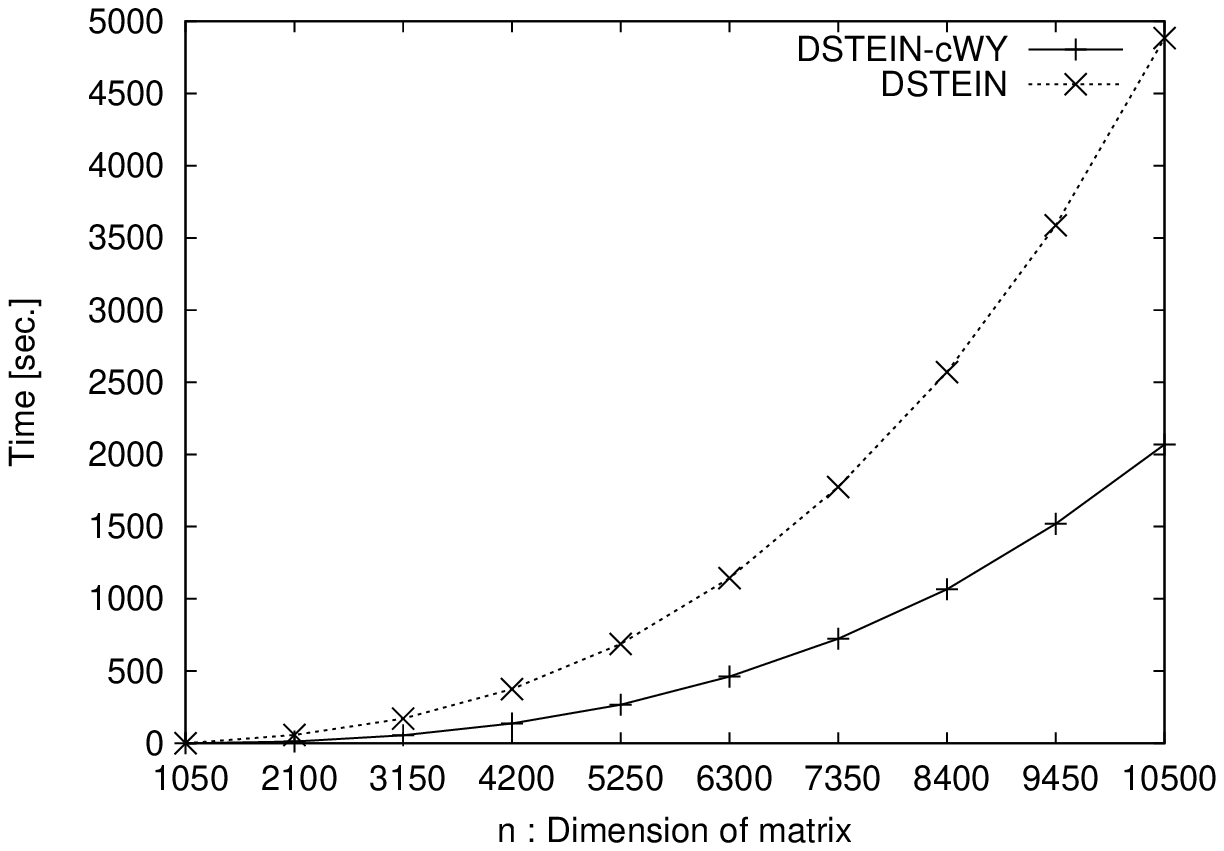}
  \end{center}
 \end{minipage}
  \caption{Dimension $n$ of Type 2 matrix and the computation time by DSTEIN and DSTEIN-cWY. 
the left graph corresponds to Computer 1 and the right Computer 2.}
  \label{graph:2}
 \begin{minipage}{0.5\hsize}
  \begin{center}
   \includegraphics[scale=.5]{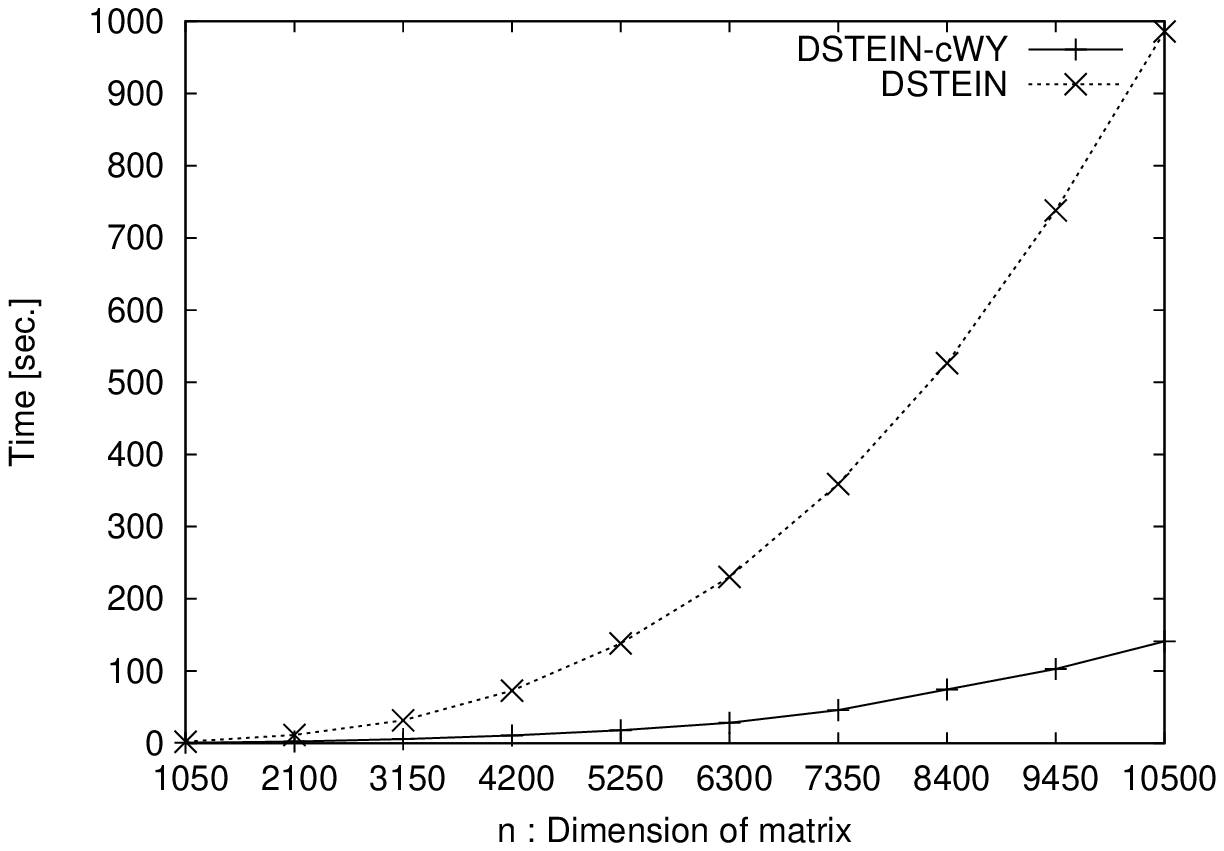}
  \end{center}
 \end{minipage}
 \begin{minipage}{0.5\hsize}
  \begin{center}
   \includegraphics[scale=.5]{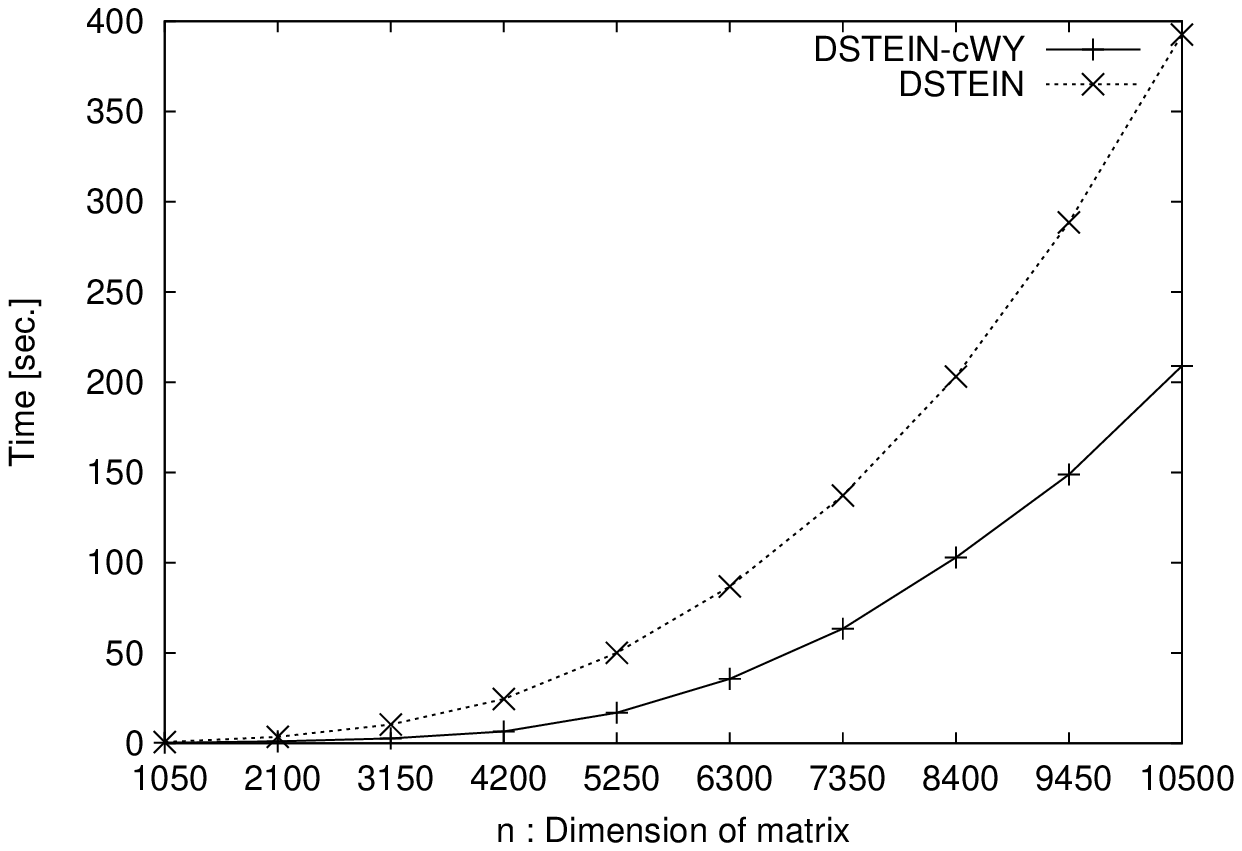}
  \end{center}
 \end{minipage}
  \caption{Dimension $n$ of Type 3 matrix and the computation time by DSTEIN and DSTEIN-cWY. 
the left graph corresponds to Computer 1 and the right Computer 2.}
  \label{graph:3}
\end{figure*}

\begin{table*}[!t]
  \footnotesize
  \centering
  \caption{Numerical results of DSTEIN and DSTEIN-cWY on Computer 1 (Type 1). }
  \label{table:3}
  \begin{tabular}{ccccccccccc} \hline
$n$ & 1050 & 2100 & 3150 & 4200 & 5250 & 6300 & 7350 & 8400 & 9450 & 10500	\\\hline\hline
$t$ [sec.] & 0.39 & 1.76 & 5.30 & 17.4 & 53.6 & 157 & 996 & 2436 & 4004 & 13231	\\
$t_{\mathrm{cwy}}$ [sec.] & 0.41 & 1.60 & 3.77 & 7.85 & 13.7 & 25.1 & 115 & 307 & 449 & 1291	\\\hline
$t/t_{\mathrm{cwy}}$        & 0.94 & 1.10 & 1.41 & 2.22 & 3.90 & 6.22 & 8.64 & 7.93 & 8.93 & 10.25	\\\hline
  \end{tabular}
  \vspace{1.0em}
  \caption{Numerical results of DSTEIN and DSTEIN-cWY on Computer 2 (Type 1). }
  \label{table:4}
  \begin{tabular}{ccccccccccc} \hline
$n$ & 1050 & 2100 & 3150 & 4200 & 5250 & 6300 & 7350 & 8400 & 9450 & 10500	\\\hline\hline
$t$ [sec.] & 0.16 & 0.75 & 2.13 & 6.41 & 19.2 & 58.3 & 372 & 889 & 1416 & 4357	\\
$t_{\mathrm{cwy}}$ [sec.] & 0.18 & 0.73 & 1.70 & 3.42 & 7.66 & 24.7 & 179 & 430 & 703 & 1933	\\\hline
$t/t_{\mathrm{cwy}}$        & 0.91 & 1.02 & 1.25 & 1.87 & 2.51 & 2.36 & 2.08 & 2.06 & 2.01 & 2.25	\\\hline
  \end{tabular}
  \vspace{1.0em}
  \caption{Numerical results of DSTEIN and DSTEIN-cWY on Computer 1 (Type 2). }
  \label{table:5}
  \begin{tabular}{ccccccccccc} \hline
$n$ & 1050 & 2100 & 3150 & 4200 & 5250 & 6300 & 7350 & 8400 & 9450 & 10500	\\\hline\hline
$t$ [sec.] & 1.73 & 154 & 448 & 989 & 1897 & 3281 & 5192 & 7749 & 10986 & 14867	\\
$t_{\mathrm{cwy}}$ [sec.] & 0.45 & 7.04 & 28.1 & 94.6 & 167 & 311 & 476 & 795 & 1029 & 1389	\\\hline
$t/t_{\mathrm{cwy}}$        & 3.85 & 21.93 & 15.94 & 10.45 & 11.34 & 10.56 & 10.92 & 9.74 & 10.68 & 10.70	\\\hline
  \end{tabular}
  \vspace{1.0em}
  \caption{Numerical results of DSTEIN and DSTEIN-cWY on Computer 2 (Type 2). }
  \label{table:6}
  \begin{tabular}{ccccccccccc} \hline
$n$ & 1050 & 2100 & 3150 & 4200 & 5250 & 6300 & 7350 & 8400 & 9450 & 10500	\\\hline\hline
$t$ [sec.] & 0.52 & 57.4 & 171 & 375 & 688 & 1143 & 1774 & 2570 & 3586 & 4884	\\
$t_{\mathrm{cwy}}$ [sec.] & 0.20 & 12.2 & 55.3 & 136 & 266 & 462 & 723 & 1067 & 1519 & 2070	\\\hline
$t/t_{\mathrm{cwy}}$        & 2.67 & 4.69 & 3.10 & 2.75 & 2.58 & 2.48 & 2.45 & 2.41 & 2.36 & 2.36	\\\hline
  \end{tabular}
  \vspace{1.0em}
  \caption{Numerical results of DSTEIN and DSTEIN-cWY on Computer 1 (Type 3). }
  \label{table:7}
  \begin{tabular}{ccccccccccc} \hline
$n$ & 1050 & 2100 & 3150 & 4200 & 5250 & 6300 & 7350 & 8400 & 9450 & 10500	\\\hline\hline
$t$ [sec.] & 2.26 & 11.5 & 31.8 & 72.9 & 138 & 230 & 359 & 526 & 738 & 986	\\
$t_{\mathrm{cwy}}$ [sec.] & 0.62 & 2.49 & 5.82 & 10.9 & 18.1 & 28.4 & 45.9 & 74.5 & 103 & 141	\\\hline
$t/t_{\mathrm{cwy}}$        & 3.66 & 4.62 & 5.47 & 6.71 & 7.66 & 8.10 & 7.82 & 7.06 & 7.18 & 6.99	\\\hline
  \end{tabular}
  \vspace{1.0em}
  \caption{Numerical results of DSTEIN and DSTEIN-cWY on Computer 2 (Type 3). }
  \label{table:8}
  \begin{tabular}{ccccccccccc} \hline
$n$ & 1050 & 2100 & 3150 & 4200 & 5250 & 6300 & 7350 & 8400 & 9450 & 10500	\\\hline\hline
$t$ [sec.] & 0.68 & 3.58 & 10.4 & 24.5 & 50.1 & 86.8 & 137 & 203 & 289 & 393	\\
$t_{\mathrm{cwy}}$ [sec.] & 0.27 & 1.10 & 2.72 & 6.59 & 16.9 & 35.7 & 63.4 & 103 & 149 & 209	\\\hline
$t/t_{\mathrm{cwy}}$        & 2.54 & 3.27 & 3.83 & 3.72 & 2.97 & 2.43 & 2.16 & 1.97 & 1.94 & 1.88	\\\hline
  \end{tabular}
\end{table*}
Table \ref{table:3}-\ref{table:8} show 
the results of the experiments on Computer 1 and 2 that are mentioned in the previous section, 
In tables, $n$ is the dimension of the experimental matrices, 
$t$ and $t_{\mathrm{cwy}}$ are computation time  by DSTEIN and DSTEIN-cWY, respectively. 
In addition, 
Fig. \ref{graph:1}-\ref{graph:3} illustrate the results 
in Tables \ref{table:3} and \ref{table:4}, \ref{table:5} and \ref{table:6}, 
\ref{table:7} and \ref{table:8} through graphs, respectively. 
In Fig. \ref{graph:1}-\ref{graph:3}, 
the dotted line corresponds to $t$ and the straight line to $t_{\mathrm{cwy}}$.\par
It is noted that DSTEIN-cWY is faster than DSTEIN 
for any cases of the all types matrices, 
without the cases of Type 1 matrix for $n=1050$. 
We see that the change from MGS to the compact WY orthogonalization 
on the DSTEIN code in parallel computing results in a significant reduction of computation time. 
We introduce a barometer $t/t_{\mathrm{cwy}}$ of the reduction effect 
by using the program DSTEIN-cWY which depends on $n$, the dimension of the experimental matrix. 
On Computer 1, 
the maximum value of $\alpha=t/t_{\mathrm{cwy}}$ is 
$\alpha=10.25$ for $n=10,500$ of Type 1,  
$\alpha=10.92$ for $n=7,350$ of Type 2,  
and $\alpha=8.10$ for $n=6,300$ of Type 3. 
On Computer 2, 
$\alpha=2.51$ for $n=5,250$ of Type 1,  
$\alpha=4.69$ for $n=2,100$ of Type 2,  
and $\alpha=3.83$ for $n=3,150$ of Type 3. 
Considering these facts, 
even if the dimension of the experimental matrices is larger than that in these examples, 
we cannot expect that the computation time can be further shortened by using DSTEIN-cWY. 
\subsection{Discussion on numerical experiments}
It is shown that DSTEIN-cWY is faster than DSTEIN 
for any dimension $n$ of the experimental matrix both on Computers 1 and 2. 
As mentioned earlier, 
according to the theoretical background in Section 3.3, 
this result shows that the compact WY orthogonalization 
is an effective algorithm for parallel computing. \par
The cause of this is related to the time required for floating-point arithmetic 
and for synchronization in parallel computing. 
The floating-point computation time increases with increasing the dimension $n$ of matrices. 
In comparison, 
the synchronization cost does not change significantly 
even if $n$ becomes larger. 
Therefore, in parallel computing, 
DSTEIN, which contains MGS (for which the number of synchronizations is large), 
creates a huge bottleneck for the synchronization cost when $n$ is small. 
This bottleneck gradually becomes less when $n$ is larger. 
However, DSTEIN-cWY has a smaller bottleneck for the synchronization cost 
because the compact WY orthogonalization requires less synchronization, 
and the floating-point computation time becomes greater than that of DSTEIN. 
This reduction effect can be seen in Table \ref{table:3}-\ref{table:8}. 
\section{Conclusions}
In this study, 
we present a new inverse iteration algorithm 
for computing all the eigenvectors of a real symmetric tri-diagonal matrix. 
The new algorithm is equipped with 
the new implementation of the compact WY orthogonalization algorithm, 
established in this paper, in the orthogonalization process. \par
Now we use a new implementation of the compact WY orthogonalization. 
Introducing this implementation, 
the computational cost of the compact WY orthogonalization can be reduced. \par
We have given numerical experiments 
for computing  eigenvectors of certain real symmetric tri-diagonal matrices 
that have many clusters with several thousand dimensions 
by using two types of inverse iteration algorithms on parallel computers. 
The results show that the compact WY inverse iteration 
is more efficient than the classical one owing to the reduction in computation time 
because of the parallelization efficiency. 
As the number of cores of the CPU increases, the parallelization efficiency increases. \par
It may be expected to apply the new inverse iteration algorithms 
to other types of matrix eigenvector problem, 
such as eigenvectors of a real symmetric band matrix, 
or singular vectors of a bidiagonal matrix. 
%

\end{document}